\documentclass[reqno,a4paper,12pt]{amsart} 

\usepackage{latexsym}

\def\serieslogo@{}

\def\@setcopyright{}

\newcommand{\N}{{\mathbb N}}
\newcommand{\R}{{\mathbb R}}
\newcommand{\sphere}{{\mathbb S}}
\newcommand{\Vol}{{\operatorname{Vol}}}
\newcommand{\dist}{{\operatorname{dist}}}

\numberwithin{equation}{section}
\theoremstyle{plain}
\newtheorem{thm}{Theorem}[section]
\newtheorem{prop}[thm]{Proposition}
\newtheorem{lemma}[thm]{Lemma}

\theoremstyle{definition}
\newtheorem{defn}[thm]{Definition}
\newtheorem{remark}[thm]{Remark}

\newenvironment{pf*}[1]{\proof[#1]}{\endproof}

\title[Wavefunctions and Densities]{Electron Wavefunctions 
       and Densities for Atoms}
\author[M. \& T. Hoffmann-Ostenhof and 
       T. \O. S\o rensen]{Maria Hoffmann-Ostenhof, 
       Thomas Hoffmann-Ostenhof and
       Thomas \O stergaard S\o rensen}

\address[M. Hoffmann-Ostenhof]{Institut f\"ur Mathematik            \\
         Strudlhofgasse 4                                           \\
         Universit\"at Wien                                         \\
         A-1090 Vienna                                              \\
         Austria}
\address[T. Hoffmann-Ostenhof]{Institut f\"ur Theoretische Chemie   \\
         W\"ahringer\-strasse 17                                    \\
         Universit\"at Wien                                         \\
         A-1090 Vienna                                              \\
         Austria} 
\address[T. Hoffmann-Ostenhof, T. \O. S\o rensen]{The Erwin 
         Schr\"odinger International 
         Institute for Mathematical Physics                         \\
         Boltzmanngasse 9                                           \\
         A-1090 Vienna                                              \\
         Austria}  

\email[M. Hoffmann-Ostenhof]{mhoffman@esi.ac.at}
\email[T. Hoffmann-Ostenhof]{thoffman@esi.ac.at}
\email[T. \O. S\o rensen]{tsorense@esi.ac.at}

\date{\today}

\thanks{Work supported by Ministerium f\"ur Wissenschaft und
  Verkehr der Republik \"Osterreich, the Austrian Science Foundation,
  grantnumber P12865-MAT, and by the European Union TMR grant FMRX-CT
  96-0001 }

\begin{document}

\begin{abstract} 
  With a special `Ansatz' we analyse the regularity pro\-perties of
  atomic electron wavefunctions and electron densities. 
  In particular we prove an a priori estimate, 
  $\sup_{y\in B(x,R)}|\nabla\psi(y)| 
   \leq  C(R)\,\sup_{y\in B(x,2R)}|\psi(y)|$ 
  and obtain for the spherically averaged electron density, 
  $\widetilde\rho(r)$, that  $\widetilde\rho''(0)$ exists and is non-negative.
\end{abstract}

\maketitle

\noindent
{\it
{Avec un `Ansatz' sp\'ecial nous analysons les propriet\'ees de 
r\'egularit\'e 
des fonctions d'ondes atomiques et des densit\'es d'electron. 
En parti\-culier nous prouvons une \'estimation a priori, 
$$\sup_{y\in B(x,R)}|\nabla\psi(y)| 
   \leq  C(R)\,\sup_{y\in B(x,2R)}|\psi(y)|$$
et obtient pour la densit\'e d'electron moyenn\'ee sur la sphere 
$\widetilde\rho(r)$, que  $\widetilde\rho''(0)$ existe et est non-n\'egative.
}}

\section{Introduction and Results}

Let \(V\) be the Coulomb potential for an atom consisting of a nucleus 
of charge \(Z\) (fixed at the origin) and \(N\) electrons:
\begin{gather}
  \label{eq:potential}
  V({\bf  x})= V(x_{1},\ldots,x_{N})=\sum_{j=1}^{N}{}-\frac{Z}{|x_{j}|}
  +\sum_{1\leq j<k\leq N}\frac{1}{|x_{j}-x_{k}|},
  \nonumber\\
  {\bf x}=(x_{1},\ldots,x_{N})\in\R^{3N},
  x_{j}=(x_{j,1},x_{j,2},x_{j,3})\in\R^{3}, j=1,\ldots,N,
\end{gather}
and let \(H\) be the corresponding \(N\) - electron Hamilton operator:
\begin{align}
  \label{eq:def-H}
  H\equiv H^{N}={}-\Delta+V
\end{align}
with
\begin{align}
  -\Delta=\sum_{j=1}^{N}-\Delta_{j}\quad,\quad
  \Delta_{j}=\sum_{i=1}^{3}\frac{\partial^{2}}{\partial x_{j,i}^{2}}
  \nonumber
\end{align}
being the kinetic energy operator of the \(N\) electrons. The quadratic
form domain of \(H\) is \(W^{1,2}(\R^{3N})\), see Reed and Simon~\cite{R&S2}. 
Assume \(\psi\in L^{2}(\R^{3N})\) is a real-valued normalised
eigenfunction of the operator \(H\):
\begin{align}
  \label{eq:eigen}
  (H-E)\psi=0\quad,\quad\|\psi\|\equiv\|\psi\|_{L^{2}(\R^{3N})}=1.
\end{align}
It is known that then \(\psi\) is continuous with bounded derivatives, 
and \(\psi\in W^{2,2}(\R^{3N})\) (Kato \cite{MR19:501a}) and that 
\(\psi\) is in fact analytic away from the singularities 
(in \(\R^{3N}\)) of \(V\), since \(V\) is here real analytic 
(see Hopf~\cite{Hopf}). In this paper we derive various qualitative 
and quantitative properties of the wave function \(\psi\) and of 
the corresponding one-electron density
\begin{align}
  \label{eq:def-rho}
  \rho(x)=\int_{\R^{3(N-1)}}|\psi(x,x_{2},\ldots,x_{N})|^{2}\,
  dx_{2}\ldots dx_{N}\quad,\quad x\in\R^{3},
\end{align}
as well as of its spherical average 
(\(x=r\omega, r=|x|, \omega=x/|x|\in\mathbb{S}^{2}\))
\begin{align}
  \label{eq:def-rho-tilde}
  \widetilde\rho(r)&=\int_{\mathbb{S}^{2}}\rho(r\omega)\,d\omega
  \nonumber\\           
  &=\int_{\mathbb{S}^{2}}\int_{\R^{3(N-1)}}
  |\psi(r\omega,x_{2},\ldots,x_{N})|^{2}
  \,dx_{2}\cdots\,dx_{N}\,d\omega \quad,\quad r\in[0,\infty).
\end{align}
\begin{remark}
\label{importance_rho}
Aside from Kato's classical results (see Kato~\cite{MR19:501a}), the local
behaviour of electron wavefunctions has been investigated more
recently by Hoffmann-Ostenhof et al.~\cite{MR95c:81171},
\cite{MHO-Seiler}. 
The electron density itself has been studied extensively in the
large-\(Z\)-limit, see Lieb~\cite{Li-Selecta}. Except for the spatial
asymptotics, see 
Ahlrichs et al.~\cite{MR82j:81064}, there are virtually no recent
rigorous results on \(\rho\) despite the fact that the density is the
central object in various popular numerical approximation schemes, as
Density Functional Theory (DFT) and all the various descendants of
Hartree-Fock theory.
\end{remark}

We now present our results.

\begin{thm}
  \label{thm:grad_est}
  Let \(\psi\) be as in \eqref{eq:eigen}. For all \(R\in(0,\infty)\),
  there exists a constant \(C=C(R)\) such that
  \begin{align}
    \sup_{{\bf y}\in B({\bf x},R)}|\nabla\psi({\bf y})|\leq C\!\!\!
    \sup_{{\bf y}\in B({\bf x},2R)}
    |\psi({\bf y})| \quad\text{ for all } \,{\bf x}\in\R^{3N}.
    \nonumber
  \end{align}
\end{thm}
\begin{remark}
  This result complements the result by 
  Simon~\cite[Thm. C.2.5 (C14)]{Si-semi} for the case of operators 
  of the form \eqref{eq:def-H}, but with \(V\) in the Kato-class
  \(K^{n,1}(\R^{n})\): for \(\delta\in[0,2)\) \((\delta=0:n\geq3)\), 
  \begin{align}
    \nonumber
    V\in K^{n,\delta}(\R^{n})\Leftrightarrow
    \lim_{\epsilon\to0}\sup_{x\in\R^{n}}
    \int_{|x-y|<\epsilon}\frac{|V(y)|}{|x-y|^{n-2+\delta}}\,dy=0.
  \end{align}
  The Coulomb potential (1.\ref{eq:potential}) is in 
  \(K^{3N,\delta}(\R^{3N})\) for all \(\delta\in[0,1)\), 
  but is not in \(K^{3N,1}(\R^{3N})\). 
\end{remark}

We recall the definition of H\"older continuity:
\begin{defn}
  \label{def:Holder}
  For \(\Omega\subset\R^{n}\) an open set, \(k\in\N\), and 
  \(\alpha\in(0,1]\), we say that the function \(u\) belongs to
  \(C^{k,\alpha}_{\text{loc}}(\Omega)\) whenever
  \(u\in C^{k}(\Omega)\), and for all \(\beta\in\N^{n}\) with 
  \(|\beta|=k\), and all open balls \(B(x_{0},r)\subset\Omega\), we have
  \begin{align*}
    \sup_{x,y\in B(x_{0},r),\,x\neq y}
    \!\!\!\!\!\!\!\!\!
    \frac{|D^{\beta}u(x)-D^{\beta}u(y)|}{|x-y|^{\alpha}} 
    \leq C(x_{0},r). 
  \end{align*}
\end{defn}

As a consequence of the proof of Theorem~\ref{thm:grad_est} we get:
\begin{prop}
  \label{prop:regularity}
  Let
  \begin{align}
     \label{eq:def_F}
     F({\bf x})=F(x_{1},\ldots,x_{N})=\sum_{j=1}^{N}{}-\frac{Z}{2}|x_{j}|
     +\sum_{1\leq j<k\leq N}\frac{1}{4}|x_{j}-x_{k}|.
  \end{align}
  Then the eigenfunction \(\psi\) given in \eqref{eq:eigen} 
  can be represented as 
  \begin{align}
    \nonumber
    \psi=e^{F}\phi
  \end{align}
  with
  \begin{align}
    \phi\in C^{1,\alpha}_{\text{{\rm loc}}}(\R^{3N})
    \text{ for all } \alpha\in(0,1).
    \nonumber
  \end{align}
\end{prop}
\begin{remark}
  This result classifies the singularities of \(\nabla\psi\) as those
  coming from \(\nabla F\): \(\nabla\psi=\psi\nabla F+e^{F}\nabla\phi\).
  Kato~\cite{MR19:501a} proved that \(\nabla\psi\) is bounded, 
  but as the ground state of Hydrogen-like systems (\(N=1, E=-Z^{2}/4\), 
  \(\psi(x)=c_{0}e^{-Z|x|/2}, x\in\R^{3}\)) shows, it is not in 
  general continuous.
\end{remark}
\begin{remark}
  \label{rem:molecules}
  The results of Theorem~\ref{thm:grad_est} 
  and Proposition~\ref{prop:regularity} easily generalise to the case
  of molecules: \(L\) nuclei, of charge \(Z_{l}\), 
 fixed at \(R_{l}\in\R^{3}\),
  \(l=1,\ldots,L\), with corresponding \(N\) - electron Hamilton operator
  \begin{align}
    H^{N,L}=\sum_{j=1}^{N}
    \left(
      {}-\Delta_{j}-\sum_{l=1}^{L}\frac{Z_{l}}{|x_{j}-R_{l}|}
    \right)
    +\sum_{1\leq j<k\leq N}\frac{1}{|x_{j}-x_{k}|}.
    \nonumber
  \end{align}
\end{remark}

  We assume throughout when studying \(\rho\) and \(\widetilde\rho\) that
  \(E\) and \(\psi\) in \eqref{eq:eigen} are such that
  there exist constants \(C_{0},\gamma>0\) such that
  \begin{align}
    \label{eq:exp_decay}
    |\psi({\bf x})|\leq C_{0}\,e^{-\gamma|{\bf x}|} \text{ for all }
    {\bf x}\in\R^{3N}.
  \end{align}
  For references on the exponential decay of eigenfunctions, see 
  e.\ g. Simon~\cite{Si-semi}.

\begin{remark}
  \label{rem:1}
  Inequality \eqref{eq:exp_decay} holds 
  when \(E<\inf \sigma_{\text{ess}}(H^{N})\).
  In this case, we let \(\varepsilon\equiv E_{0}^{N-1}-E\) with
  \(E_{0}^{N-1}\) the ground state energy of the \((N-1)\) - electron
  operator:
  \begin{align}
    \label{def:H_N-1}
    H^{N-1}= \sum_{j=2}^{N}\left(-\Delta_{j}
    -\frac{Z}{|x_{j}|}\right)
    +\sum_{2\leq j<k\leq N}\frac{1}{|x_{j}-x_{k}|}.
  \end{align}
  By the HVZ-theorem (see Cycon et al.~\cite[Theorem 3.7]{CFKS}, 
  \(\inf\sigma_{\text{ess}}(H^{N})=E^{N-1}_{0}\), 
  and so \(\varepsilon>0\) if \(E<\inf\sigma_{\text{ess}}(H^{N})\). 
  When we study \(H\) in a symmetry sector, so that \(\psi\)
  transforms according to this symmetry, then \(E_{0}^{N-1}\) stands
  for the groundstate energy of the ionized particle system described
  by the Hamiltonian \(H^{N-1}\) in the appropriate symmetry subspace
  as determined by the symmetry behaviour of \(\psi\). For this case a 
  modified version of the HVZ-theorem holds (see e.\ g. Reed
  and Simon~\cite[Thm. XIII. 17']{R&S4} and Zhislin and
  Sigalov~\cite{Zh-Si}). This includes in particular the physically
  important case of real atoms (Pauli principle). So if \(E\) lies
  below the beginning of the essential spectrum of \(H\) considered in 
  a symmetry sector, then analogously to the above the ionisation
  energy \(\varepsilon>0\) and \(\psi\) satisfies \eqref{eq:exp_decay}.
\end{remark}
\begin{remark}
  When assuming \eqref{eq:exp_decay}, Theorem~\ref{thm:grad_est} 
  implies that \(|\nabla\psi({\bf x})|\) also decays exponentially 
  for \(|{\bf x}|\to\infty\).
\end{remark}
\begin{remark}
  From \eqref{eq:exp_decay} and Lebesgue's Dominated Convergence 
  \(\,\,\)
  Theorem
  follows that the density \(\rho\) is continuous in
  \(\R^{3}\).
\end{remark}

\renewcommand{\theenumi}{{\rm \roman{enumi}}}

\begin{thm}
  \label{thm:eq-rho-tilde}
  Let \(\psi\) be given according to \eqref{eq:eigen} and assume
  that \eqref{eq:exp_decay} holds. Then:
  \begin{enumerate}
  \item
  The function \(\rho\) defined in \eqref{eq:def-rho} satisfies, in the  
  distributional sense, the equation
  \begin{align}
    \label{eq:schr-like}
    {}-\frac{1}{2}\Delta\rho-\frac{Z}{r}\rho+h=0 
    \quad\text{ in }\quad \R^{3},
  \end{align}
  where 
  \begin{align}
    \nonumber
    h&\in C^{\alpha}(\R^{3}\setminus\{0\})\cap
    L^{\infty}(\R^{3})\,\,
    \text{ for all } \,\,  \alpha\in(0,1)
    \intertext{ and }
    \rho&\in
    C^{2,\alpha}(\R^{3}\setminus\{0\})\cap C^{0,1}(\R^{3})
    \,\,\text{ for all }\,\,
    \alpha\in(0,1).
    \nonumber
  \end{align} 
  \label{No1}   
  \item
  The function \(\widetilde\rho\) defined in \eqref{eq:def-rho-tilde} 
  satisfies  
  \begin{align}
    \label{eq:schr-like-tilde}
    {}-\frac{1}{2}\Delta\widetilde\rho-\frac{Z}{r}\widetilde\rho+\widetilde h
    =0 \text{ for } r>0
  \end{align}
  where \(\widetilde
  h(r)=\int_{\sphere^{2}}h(r\omega)\,d\omega\). Thereby,
  \begin{align}
    \nonumber
    \widetilde h\in C^{\alpha}((0,\infty))\cap
    C^{0}([0,\infty))
    \text{ for all } \alpha\in(0,1)
  \end{align}
  and
  \begin{align}
    \nonumber
    \widetilde\rho\in C^{2,\alpha}((0,\infty))\cap C^{2}([0,\infty))
    \text{ for all } \alpha\in(0,1).
  \end{align}
  \item
  \begin{align}
    \label{eq:up-bd}
    h(x)&\leq C(R) 
    \left(
         \int_{B(x,R)}\!\!\!\!\!\!\!\!\rho(y)\,dy         
         +\rho(x)
    \right)
    \,\,\,\text{ for all } x\in\R^{3}, 
    \\
    \label{eq:low-bd}
    h(x)&\geq\varepsilon\,\rho(x)
    \text{ for all } x\in\R^{3} , \text { if }
    \varepsilon=E_{0}^{N-1}-E>0.
  \end{align}
  \item
  \begin{align}
    \label{eq:formel-rho-tilde}
    \Big(\frac{d^{2}}{dr^{2}}\,\widetilde\rho\Big)(0)=\frac{2}{3}\big(\widetilde
    h(0)+Z^{2}\widetilde\rho(0)\big).
  \end{align}
  \end{enumerate}
\end{thm}

\ %

\begin{remark}
  The results in (i) generalize to the case of
  molecules, where the continuity results for \(\rho\) and \(h\) hold
  in the complement of the set \(\{R_{1},\ldots,R_{L}\}\subset\R^{3}\) (see Remark~\ref{rem:molecules}).
\end{remark}

\begin{remark}
It is known that eigenfunctions obey (Kato's) Cusp Condition (see
Kato~\cite{MR19:501a}), and similar properties hold for particle
densities. For more recent results see Hoffmann-Ostenhof et 
al.~\cite{MR95c:81171},~\cite{MR58:25667}, Hoffmann-Ostenhof and
Seiler~\cite{MHO-Seiler}. In the proof of 
Theorem~\ref{thm:eq-rho-tilde}, (iv) we make use of the Cusp
Condition for \(\widetilde\rho\), namely:
\begin{align}
  \widetilde\rho'(0)=\lim_{r\downarrow0}\frac{\widetilde\rho(r)-\widetilde\rho(0)}{r}
  =-Z\widetilde\rho(0) \text{ and }
  \lim_{r\downarrow0}\widetilde\rho'(r)=\widetilde\rho'(0)
  \label{eq:Cusp}
\end{align}
and also present a proof for it.
\end{remark}

\begin{remark}
Of course our results are only first steps in a thorough
investigation of qualitative properties of the one-electron
density. Here are some obvious open questions:
\begin{enumerate}
   
\item
Is \(\rho(x)>0\) for all \(x\in\R^{3}\)? We remark that this
cannot be true in general, since it is false for some exited states of 
Hydrogen.
 
\item
Is \(\rho\in C^{\infty}(\R^{3}\setminus\{0\})\) or even
\(C^{\omega}(\R^{3}\setminus\{0\})\) ?

\item
Is \(\widetilde\rho\) smooth in \([0,\infty)\), in the sense that
\(\big(\frac{d^{k}}{dr^{k}}\widetilde\rho\big)(r)\) exists for \(r\geq0\) 
for all \(k\)? 

\item
Is \(\frac{d}{dr}\widetilde\rho(r)\leq0\) for \(r\geq0\) ?
This is expected to be true for groundstate densities, but not known
even for the bosonic case like Helium. Our results imply that 
\(\frac{d}{dr}\widetilde\rho(r)\leq0\) for \(r\leq R_{0}\) for the bosonic 
case, where \(R_{0}\) depends on the constant \(C\) in
Theorem~\ref{thm:grad_est}. Note that because of 
\eqref{eq:schr-like} and \eqref{eq:low-bd} we have
\(\Delta\rho\geq0\) for \(|x|\geq Z/\varepsilon\), and so the Maximum
Principle gives that  \(\frac{d}{dr}\widetilde\rho(r)<0\) for \(r>
Z/\varepsilon\).  
\end{enumerate}
\end{remark}
\begin{remark}
\label{t_1}
In the proof of Theorem~\ref{thm:eq-rho-tilde} we obtain (see
Proposition~\ref{prop:reg}):
With \(\nabla_{1}=(\frac{\partial}{\partial x_{1,1}},
\frac{\partial}{\partial x_{1,2}},\frac{\partial}{\partial x_{1,3}})\),
the function
\begin{align}
  \nonumber
  t_{1}(r)=\int_{\mathbb{S}^{2}}
  \int_{\R^{3(N-1)}}|\nabla_{1}\psi(r\omega,x_{2},\ldots,x_{N})|^{2}\,
  dx_{2}\ldots dx_{N}\,d\omega
  \end{align}
  is continuous on \([0,\infty)\).
\end{remark}

\section{Proofs}

Throughout the proofs, we will denote by \(C\) generic constants.

Crucial for our investigations is Corollary~8.36 in Gilbarg and
Tru\-dinger~\cite{GandT}. We shall make use of this 
result several times and for convenience we state it already here, 
adapted for our special case:
\begin{prop}
  \label{prop:GT}
  Let \(\Omega\) be a bounded domain in \(\R^{n}\) and suppose \(u\in
  W^{1,2}(\Omega)\) is a weak solution of \(\Delta u
  +\sum_{j=1}^{n}b_{j}D_{j}u+Wu=g\) in \(\Omega\), where
  \(b_{j},W,g\in L^{\infty}(\Omega)\). Then \(u\in
  C^{1,\alpha}(\Omega)\) for all \(\alpha\in(0,1)\) and for any domain
  \(\Omega'\), \(\overline{\Omega'}\subset\Omega\) we have
  \begin{align}
    |u|_{C^{1,\alpha}(\Omega')}\leq
    C\big(\sup_{\Omega}|u|+\sup_{\Omega}|g|\big) 
    \nonumber
  \end{align}
  for \(C=C(n,M,\dist(\Omega',\partial\Omega))\), with 
  \begin{align}
    \max_{j=1,\ldots,n}\{1,\|b_{j}\|_{L^{\infty}(\Omega)},
    \|W\|_{L^{\infty}(\Omega)},\|g\|_{L^{\infty}(\Omega)}\} \leq M.
    \nonumber
  \end{align}
  Thereby 
  \begin{align}
    \nonumber
    |u|_{C^{1,\alpha}(\Omega')}=\|u\|_{L^{\infty}(\Omega')}
    +\|\nabla u\|_{L^{\infty}(\Omega')}
    +\sup_{x,y\in\Omega',\,x\neq y}
    \frac{|\nabla u(x)-\nabla u(y)|}{|x-y|^{\alpha}}.
  \end{align}
\end{prop}

\begin{center}
  Proof of Theorem~\ref{thm:grad_est} and Proposition~\ref{prop:regularity}.
\end{center}

\ %

Let the function \(F\) be as in \eqref{eq:def_F} and define the function 
\(F_{1}\) by
\begin{align}
  \label{def_F_1}  
  F_{1}(x_{1},\ldots,x_{N})
  =\sum_{j=1}^{N}{}-\frac{Z}{2}\sqrt{|x_{j}|^{2}+1}
  +\sum_{1\leq j<k\leq N}\frac{1}{4}
  \sqrt{|x_{j}-x_{k}|^{2}+1}. 
\end{align}
A computation shows that
\begin{align}
  \label{eq:delta_F}
  &\Delta F=V,
  \\
  \label{eq:est_F}
  &\|F-F_{1}\|_{L^{\infty}(\R^{3N})}
  \,\,,\,\,\,
  \|\nabla(F-F_{1})\|_{L^{\infty}(\R^{3N})}
  \,\,,\,\,\,
  \|D^{\beta} F_{1}\|_{L^{\infty}(\R^{3N})}
  \nonumber\\
  &\leq
  C(\beta,N,Z)\,\,,\,\,\, |\beta|\geq1.
\end{align}
Make the `Ansatz' \(\psi=e^{F-F_{1}}\psi_{1}\). Using
\((H-E)\psi=0\) and \eqref{eq:delta_F} we get that \(\psi_{1}\)
satisfies the equation
\begin{align}
  \label{eq:elliptic_eq}
  \Delta\psi_{1}+2\nabla(F-F_{1})\cdot\nabla\psi_{1}
  +(|\nabla(F-F_{1})|^{2}-\Delta F_{1}+E)\psi_{1}=0.
\end{align}
Due to \eqref{eq:est_F} the coefficients in \eqref{eq:elliptic_eq}
are bounded in \(\R^{3N}\). Then Proposition~\ref{prop:GT} implies that 
\(\psi_{1}\) is \(C^{1,\alpha}\) for all
\(\alpha\in(0,1)\), in any ball \(B({\bf x},R)\subset\R^{3N}\),
and
\begin{align}
  \label{eq:apriori_1}
  |\psi_{1}|_{C^{1,\alpha}(B({\bf x},R))}
  \leq C\, \sup_{y\in B({\bf x},2R)}|\psi_{1}(y)|
\end{align}
with \(C\) depending on \(R\) but not on \({\bf x}\).
Since
\begin{align}
  &|\nabla\psi({\bf y})|
  \leq
  |\nabla(F-F_{1})|\,|\psi({\bf y})| +
  |e^{F-F_{1}}\nabla\psi_{1}({\bf y})|
  \nonumber
\end{align}
we obtain, via \eqref{eq:est_F} and \eqref{eq:apriori_1},
\begin{align}
  &{}
  \sup_{{\bf y}\in B({\bf x},R)}|\nabla\psi({\bf y})|
  \leq C\big(\!\!\!\!
  \sup_{{\bf y}\in B({\bf x},R)}|\psi({\bf y})|
  +\sup_{{\bf y}\in B({\bf x},R)}|\nabla\psi_{1}({\bf y})|\big)
  \nonumber\\
  &\leq C\big(\!\!\!\!
  \sup_{{\bf y}\in B({\bf x},2R)}|\psi({\bf y})|
  +\sup_{{\bf y}\in B({\bf x},2R)}|\psi_{1}({\bf y})|
  \big)\leq C\!\!\!
  \sup_{{\bf y}\in B({\bf x},2R)}|\psi({\bf y})|,
  \nonumber
\end{align}
with \(C=C(R)\).
This proves Theorem~\ref{thm:grad_est}.

Proposition~\ref{prop:regularity} follows from
\(\psi=e^{F-F_{1}}\psi_{1}\), 
\(\psi_{1}\in C^{1,\alpha}_{\text{{\rm loc}}}(B({\bf x},R))\), since 
\(e^{-F_{1}}\) is smooth.
\qed

\ %

\begin{center}
  Proof of Theorem~\ref{thm:eq-rho-tilde}.
\end{center}
  Multiplying the equation \((H-E)\psi=0\) with \(\psi\) 
  and integrating over \(x_{2},\ldots,x_{N}\) gives the equation
  \begin{align}
    \label{eq:H_N-1}
    \int_{\R^{3(N-1)}}&\psi\Delta_{1}\psi\,dx_{2}\cdots\,dx_{N}
    +\frac{Z}{|x_{1}|}\,\rho(x_{1})=\nonumber\\&=\sum_{j=2}^{N}
    \int_{\R^{3(N-1)}}\psi\big(-\Delta_{j}-\frac{Z}{|x_{j}|}\big)
    \psi\,dx_{2}\cdots\,dx_{N}\nonumber\\&\quad
    +\sum_{1\leq j<k\leq N}\int_{\R^{3(N-1)}}
    \frac{1}{|x_{j}-x_{k}|}\,\psi^{2}\,dx_{2}\cdots\,dx_{N}
    \nonumber\\&=\int_{\R^{3(N-1)}}\psi\big(H^{N-1}-E\big)\psi
    \,dx_{2}\cdots\,dx_{N}\nonumber\\&\quad
    +\sum_{j=2}^{N}\int_{\R^{3(N-1)}}\frac{1}{|x_{1}-x_{j}|}\,\psi^{2}
    \,dx_{2}\cdots\,dx_{N}
  \end{align}
  where \(H^{N-1}\) is the \((N-1)\) - electron operator defined in
  \eqref{def:H_N-1}.
  Since \(\Delta_{1}\big(\psi^{2}\big)=2|\nabla_{1}\psi|^{2}
  +2\psi\Delta_{1}\psi\)
  and \(\int\Delta_{1}(\psi^{2})(x_{1},x')\,dx'=\Delta_{1}\rho\) in
  the distributional sense, we get that
  \begin{align}
    \frac{1}{2}\Delta_{1}\rho(x_{1})&=\int_{\R^{3(N-1)}}
    |\nabla_{1}\psi|^{2}\,dx_{2}\cdots\,dx_{N}\nonumber\\
    \label{eq:delta_rho}&\quad+\int_{\R^{3(N-1)}}
    \psi\Delta_{1}\psi\,dx_{2}\cdots\,dx_{N}
  \end{align}
  which, together with \eqref{eq:H_N-1}, gives the equation
  (\(r=|x|\))
  \begin{align}
    \label{eq:rho_first}
    \frac{1}{2}\Delta\rho(x)+\frac{Z}{r}\rho(x)
    &=\int_{\R^{3(N-1)}}\psi\big(H^{N-1}-E\big)\psi
    \,dx_{2}\cdots\,dx_{N}\nonumber\\
    &\quad+\sum_{j=2}^{N}\int_{\R^{3(N-1)}}
    \frac{1}{|x_{1}-x_{j}|}\,\psi^{2}\,dx_{2}\cdots\,dx_{N}
    \nonumber\\&\quad+\int_{\R^{3(N-1)}}|\nabla_{1}\psi|^{2}
    \,dx_{2}\cdots\,dx_{N}\equiv h(x),
  \end{align}
hence we obtain \eqref{eq:schr-like}. Integration of
\eqref{eq:schr-like} over \(\sphere^{2}\) yields
\eqref{eq:schr-like-tilde}.
 
  The proof of the regularity properties of the functions \(h\) 
  and \(\widetilde h\) are rather technical, and therefore postponed 
  to the next section.

  We now verify the regularity properties of the functions
  \(\rho\) and \(\widetilde\rho\) under the assumption that the
  regularity properties of \(h\) and \(\widetilde h\) stated in
  Theorem~\ref{thm:eq-rho-tilde} have been shown.
  Define the function \(\mu\) by the equation (\(r=|x|\))
  \begin{align}
    \label{def:mu}
    \rho(x)=e^{-Zr}\big(\rho(0)+\mu(x)\big).
  \end{align}
  Then \(\mu=e^{Zr}\rho-\rho(0)\), \(\mu(0)=0\), and 
  \eqref{eq:rho_first} implies that
  \begin{align}
    \label{eq:mu}
    \Delta\mu-2Z\,\frac{x}{r}\cdot\nabla\mu+Z^{2}\mu=2h e^{Zr}-Z^{2}\rho(0).
  \end{align}
  Since \(h\in L^{\infty}_{\text{loc}}\), all coefficients of
  \eqref{eq:mu} are \(L^{\infty}_{\text{loc}}\), and since \(\rho\in
  W^{1,2}_{\text{loc}}\), also \(\mu\in
  W^{1,2}_{\text{loc}}\). Therefore Proposition~\ref{prop:GT} leads to
  \(\mu\in C^{1,\alpha}_{\text{loc}}\),
  for all \(\alpha\in(0,1)\). Due to \eqref{def:mu},
  \(\rho\in C^{0,1}(\R^{3})\) follows.

  Now consider
  \begin{align}
    \label{eq:rho2}
    \Delta\rho={}-\frac{2Z}{r}\rho+2h\equiv g\quad\text{ in } \quad
    \R^{3}\setminus \{0\}.
  \end{align}
  Since \(h\in C^{\alpha}(\R^{3}\setminus\{0\})\),
  for all \(\alpha\in(0,1)\) and due to the above,
  \(\rho/r\in C^{\alpha}(\R^{3}\setminus\{0\})\),
  for all \(\alpha\in(0,1)\), we have
  \begin{align}
    \label{eq:reg_g}
    g\in
    C^{\alpha}(\R^{3}\setminus\{0\}) \text{ for all }\alpha\in(0,1).
  \end{align}
  From \eqref{eq:rho2} and \eqref{eq:reg_g} we obtain from regularity theory
  for the Poisson equation that \(\rho\in
  C^{2,\alpha}(\R^{3}\setminus\{0\})\) (see e.\ g.\ Gilbarg and
  Trudinger~\cite[Thm. 4.3 and 4.6]{GandT} or 
  Lieb and Loss~\cite[Thm. 10.3]{Li-Loss}).

  We proceed analogously for \(\widetilde\rho\):
  integrating \eqref{eq:rho2} over \(\sphere^{2}\), we get the equation
  \begin{align}
    \label{eq:rho_tilde2}
    \Delta\widetilde\rho={}-\frac{2Z}{r}\widetilde\rho
    +2\widetilde h\equiv \widetilde g\quad\text{ in } \quad
    \R^{3}\setminus \{0\}
  \end{align}
  with
  \begin{align}
    \label{def_h_tilde}
    \widetilde h(r)=\int_{\sphere^{2}}h(r\omega)\,d\omega.
  \end{align}
  Since the R.H.S. of \eqref{eq:rho_tilde2} is in
  \(C^{\alpha}(\R^{3}\setminus\{0\})\), we obtain
  that \(\widetilde\rho\) as a (radially symmetric) function in \(\R^{3}\)
  is \(C^{2,\alpha}\) away from the origin,
  and therefore
  \(\widetilde\rho:\R_{+}\to\R\) satisfies \(\widetilde\rho\in
  C^{2,\alpha}((0,\infty))\).

  That \(\widetilde\rho\in C^{2}([0,\infty))\) is shown in the proof of (iv).

  Next we prove (iii). To prove the bound \eqref{eq:low-bd}, 
  let \(E_{0}^{N-1}\) be the groundstate energy for the operator 
  \(H^{N-1}\). From Remark~\ref{rem:1} and the Variational Principle 
  we get that for almost all \(x_{1}\in\R^{3}\),
  \begin{align}
    \int_{\R^{3(N-1)}}\psi\big(H^{N-1}-E\big)\psi
    \,dx_{2}\cdots\,dx_{N}
    \geq(E_{0}^{N-1}-E)\rho(x_{1}),
    \nonumber
  \end{align}
  and so \(h(x)\geq \varepsilon\,\rho(x)\) with
  \(\varepsilon=E_{0}^{N-1}-E>0\).

  As for the bound \eqref{eq:up-bd}, note that due to the operator
  inequality \(-\Delta-\beta/r\geq-\beta^{2}/4\) (true in dimension
  \(3\)) and the translation invariance of \({}-\Delta\) we have, 
  for almost all \(x_{k}\in\R^{3}\) (fixed), \(k\in\{1,\ldots,N\}\),
  \(k\neq j\),
 \begin{align}
    \int_{\R^{3}}\frac{1}{|x_{j}-x_{k}|}\,|\psi|^{2}\,dx_{j}
    \leq \int_{\R^{3}}|\nabla_{j}\psi|^{2}\,dx_{j}
    +\frac{1}{4}\int_{\R^{3}}|\psi|^{2}\,dx_{j}.
    \nonumber
  \end{align}
   In this way, using \eqref{eq:delta_rho},
  \begin{align}
    \label{est_h}
    h(x)\leq C\big(\int_{\R^{3(N-1)}}|\nabla\psi|^{2}\,dx_{2}\cdots dx_{N}
    + \int_{\R^{3(N-1)}}|\psi|^{2}\,dx_{2}\cdots dx_{N}\big).
  \end{align}
  Due to Theorem~\ref{thm:grad_est} and a subsolution estimate
  (see Simon~\cite[Theorem C.1.2.]{Si-semi}) we get,
  with \({\bf x}=(x_{1},\cdots,x_{N})=(x_{1},x')\),
  \(x'\in\R^{3(N-1)}\) and \(\chi_{\Omega}\) the characteristic
  function of the set \(\Omega\):
  \begin{align}
    |\nabla\psi({\bf x})|^{2}\leq C\sup_{{\bf y}\in 
    B({\bf x},R)}|\psi({\bf  y})|^{2} 
    &\leq C\int_{{\bf  y}\in B({\bf x},2R)}|\psi({\bf  y})|^{2} 
    \,d{\bf y}
    \nonumber\\
    &=C\int_{\R^{3N}}\chi_{B({\bf x},2R)}({\bf y})|\psi({\bf  y})|^{2} 
    \,d{\bf y}.
    \nonumber
  \end{align}
  Using this, and that for fixed \({\bf x,y}\in\R^{3N}\):
  \begin{align}
    \chi_{B({\bf x},2R)}({\bf y})=\chi_{B({\bf y},2R)}({\bf x})
    =\begin{cases} 
               & 1 \quad\text{ if } |{\bf x}-{\bf y}|<2R, \\
               & 0 \quad\text{ otherwise }
     \end{cases}
    \nonumber
  \end{align}
  we get, by Fubini,
  \begin{align}
    \label{chi_trick}
    &\int_{\R^{3(N-1)}}|\nabla\psi(x_{1},x')|^{2}\,dx'
    \nonumber\\&\leq C \int_{\R^{3(N-1)}}
    \Big(\int_{\R^{3N}}\chi_{B({\bf x},2R)}({\bf y})|\psi({\bf 
    y})|^{2}\,d{\bf y}\Big)\,dx'\nonumber\\
    &=C\,\int_{\R^{3N}}|\psi({\bf y})|^{2}
    \Big(\int_{\R^{3(N-1)}}
    \chi_{B({\bf y},2R)}({\bf x})\,dx'\Big)\,d{\bf y}.
  \end{align}
  Note that with \({\bf z}={\bf x}-{\bf y}\) we have
  \begin{align}
    \chi_{B({\bf y},2R)}({\bf x})=\chi_{B({\bf y},2R)}({\bf z+y})
    =\chi_{B({\bf 0},2R)}({\bf z})
    \nonumber
  \end{align}
  and so (with \(r'=|z'|\) and \(\omega'=z'/r'\))
  \begin{align}
    \label{int_chi}
    &\int_{\R^{3(N-1)}}\!\!\!\!\chi_{B({\bf y},2R)}((x_{1},x'))
    \,dx'=\int_{\R^{3(N-1)}}\chi_{B({\bf 0},2R)}((z_{1},z'))\,dz'
    \nonumber\\&\quad=\int_{|(z_{1},z')|\leq2R}
    \!\!\!\!\!\!\!\!\!\!\!\!\!\!dz'\quad
    =\chi_{B(0,2R)}(z_{1})\int_{\sphere^{3(N-1)-1}}
    \int_{0}^{\sqrt{4R^{2}-|z_{1}|^{2}}}\!\!\!\!\!
    {r'}^{3(N-1)-1}\,dr'\,d\omega'\nonumber\\
    &\quad=C(N)\big(4R^{2}-|z_{1}|^{2})^{3(N-1)/2}\chi_{B(0,2R)}(z_{1})
    \nonumber\\
    &\quad\leq \widetilde C(N) R^{3(N-1)}\chi_{B(x_{1},2R)}(y_{1}).
  \end{align}
  From \eqref{chi_trick} and \eqref{int_chi} we get
  \begin{align}
    \label{est_nabla_psi}
    &\int_{\R^{3(N-1)}}|\nabla\psi(x_{1},x')|^{2}\,dx'
    \nonumber\\
    &\leq C(R)\int_{|x_{1}-y_{1}|<2R}\int_{\R^{3(N-1)}}|\psi(y_{1},y')|^{2}
    \,dy'\,dy_{1}= C(R)\int_{B(x_{1},2R)}\!\!\!\!\!\!\!\!\!\!\!\!\!\!
    \rho(y_{1})\,dy_{1}.
  \end{align}
  Combining \eqref{est_h} and \eqref{est_nabla_psi} proves
  \eqref{eq:up-bd}.

  We now prove (iv). 
  We first prove Kato's Cusp Condition \eqref{eq:Cusp} for the function
  \(\widetilde\rho\):
  \begin{align}
    \widetilde\rho'(0)=\lim_{r\downarrow0}\frac{\widetilde\rho(r)-\widetilde\rho(0)}{r}
    =-Z\widetilde\rho(0) \text{ and }
    \lim_{r\downarrow0}\widetilde\rho'(r)=\widetilde\rho'(0).
    \nonumber
  \end{align}

  First, define the function \(\widetilde\mu\) by the equation (see also
  \eqref{def:mu})
  \begin{align}
    \label{def_eta}
    \widetilde\rho(r)=e^{-Zr}\big(\widetilde\rho(0)+\widetilde\mu(r)\big).
  \end{align}
  Note that \(\widetilde\mu(0)=0\). Then, using \eqref{eq:schr-like-tilde},
  \(\widetilde\mu\) satisfies the equation
  \begin{align}
    \Delta\widetilde\mu-2Z\,\frac{x}{r}\cdot\nabla\widetilde\mu+Z^{2}\widetilde\mu
    =2\widetilde h e^{Zr}-Z^{2}\widetilde\rho(0),
    \nonumber
  \end{align} 
  and, since \(\widetilde h\) is continuous, Proposition~\ref{prop:GT}
  gives  that \(\widetilde\mu\), as a (radially symmetric)
  function in \(\R^{3}\), 
  is \(C^{1,\alpha}\) in a neighbourhood of the origin. In
  particular, \(\lim_{r\downarrow0}\widetilde\mu'(r)=\widetilde\mu'(0)\). 
  Since (see \eqref{def_eta})
  \begin{align}
  \label{eq:rho-tilde-der}  
    \widetilde\rho'(r)=-Z\widetilde\rho(0)+e^{-Zr}\widetilde\mu'(r)
  \end{align}
  this means that
  \begin{align}
    \label{eq:eta'(0)}
    \lim_{r\downarrow0}\widetilde\rho'(r)=-Z\widetilde\rho(0)+\lim_{r\downarrow0}\widetilde\mu'(r)
    =\widetilde\mu'(0)-Z\widetilde\rho(0).
  \end{align}
  From \eqref{def_eta} we also get that
  \begin{align}
    \frac{\widetilde\mu(r)}{r}=e^{Zr}\,\frac{\widetilde\rho(r)-\widetilde\rho(0)}{r}+ 
    \frac{e^{Zr}-1}{r}\,\widetilde\rho(0).
    \nonumber
  \end{align}
  This, together with \eqref{eq:eta'(0)} and \(\widetilde\mu(0)=0\), implies that
  \begin{align}
   \label{eq:cont_rho'}
    \lim_{r\downarrow0}\widetilde\rho'(r)&=\widetilde\mu'(0)-Z\widetilde\rho(0)
    =\lim_{r\downarrow0}\frac{\widetilde\mu(r)}{r}
    -\widetilde\rho(0)\,\lim_{r\downarrow0}\frac{e^{Zr}-1}{r} 
    \nonumber\\
    &=\lim_{r\downarrow0}e^{Zr}\frac{\widetilde\rho(r)-\widetilde\rho(0)}{r}
    =\widetilde\rho'(0),
  \end{align}
  and by \eqref{eq:rho-tilde-der} and \eqref{eq:cont_rho'},
  \(\widetilde\rho'(0)=-Z\widetilde\rho(0)\). This proves \eqref{eq:Cusp}.

Next, define the function \(\widetilde\eta\) by the equations
\begin{align}
  \label{eq:def_tilde_eta}
  \widetilde\rho(r)&=e^{-Zr}\big(\widetilde\rho(0)+\beta
  r^{2}+\widetilde\eta(r)\big),\\
  \beta&=\frac{1}{3}\big(\widetilde h(0)-\frac{Z^{2}}{2}\widetilde\rho(0)\big). 
  \nonumber
\end{align}
Then \(\widetilde\eta(0)=0\), and due to \eqref{eq:Cusp},
\(\widetilde\eta'(0)=0\). Together with \eqref{eq:schr-like} this gives
\begin{align}
  \label{deg:G}
  \Delta &\widetilde\eta(r)=\Big(\frac{1}{r^{2}}\frac{\partial}{\partial
    r}r^{2}\frac{\partial}{\partial r}\Big)\widetilde\eta
    \nonumber\\
    &=2\big[e^{Zr}\widetilde h(r)-\frac{Z^{2}}{2}\widetilde\rho(r)e^{Zr}
    -3\beta+2Z\beta r+Z\widetilde\eta'(r)\big]
    \nonumber\\
    &\equiv2G(r).
\end{align}
From the foregoing regularity properties of
\(\widetilde\rho\), in particular \eqref{eq:Cusp} and
\eqref{eq:def_tilde_eta}, we obtain that \(\widetilde\eta'\in
C^{0}([0,\infty))\). From this, together with the regularity properties of
\(\widetilde h\) shown in Section~\ref{app:reg}, we conclude 
that \(G\in C^{0}([0,\infty))\) and
\begin{align}
  \label{eq:limit_G}
  G(r)\to \widetilde h(0)-\frac{Z^{2}}{2}\widetilde\rho(0)-3\beta=0
  \quad\text{as}\quad r\downarrow 0.
\end{align}
From \eqref{eq:limit_G} and \eqref{deg:G} we get that
\begin{align}
  \widetilde\eta(r)=2\int_{0}^{r}\frac{1}{t^{2}}\int_{0}^{t}G(s)s^{2}\,ds\,dt,
   \nonumber
\end{align}
and
\begin{align}
  \frac{\widetilde\eta'(r)-\widetilde\eta'(0)}{r}
  =\frac{2}{r^{3}}\int_{0}^{r}G(s)s^{2}\,ds
    =\frac{2}{3}\frac{1}{\Vol_{\R^{3}}(B(0,r))}
  \int_{B(0,r)}\!\!\!\!\!\!G(|x|)\,d^{3}\!x,
  \nonumber
\end{align}
so that
\begin{align}
   \label{eq:tilde_eta''}
 \eta''(0)=\lim_{r\downarrow0}\frac{\eta'(r)-\eta'(0)}{r}
 =\frac{2}{3}\,G(0)=0.  
\end{align}
Then by \eqref{eq:def_tilde_eta}
\(\widetilde\rho''(0)\) exists, and
\begin{align}
  \widetilde\rho''(0)=Z^{2}\widetilde\rho(0)+2\beta
   =\frac{2}{3}\big(\widetilde h(0)+Z^{2}\widetilde\rho(0)\big).   
  \nonumber
\end{align}
This verifies \eqref{eq:formel-rho-tilde}.

Furthermore, by \eqref{deg:G}, 
\begin{align}
  \Delta \widetilde\eta(r)=\widetilde\eta''(r)+\frac{2}{r}\widetilde\eta'(r)=2G(r)
  \nonumber
\end{align}
and so
\begin{align}
  \nonumber
  \widetilde\eta''(r)=2G(r)-\frac{2}{r}\widetilde\eta'(r)
  =2G(r)-2\left(\frac{\widetilde\eta'(r)-\widetilde\eta'(0)}{r}\right)
\end{align}
since \(\widetilde\eta'(0)=0\). This implies, by \eqref{eq:tilde_eta''},
\begin{align}
  \lim_{r\downarrow0}\widetilde\eta''(r)
  =2\big(\lim_{r\downarrow0}G(r)
  -\lim_{r\downarrow0}\left(\frac{\widetilde\eta'(r)
  -\widetilde\eta'(0)}{r}\right)\big)
  =\widetilde\eta''(0)=0,
  \nonumber
\end{align}
so that due to \eqref{eq:def_tilde_eta} \(\widetilde\rho''(r)\) 
is continuous at \(r=0\). 
Hence formula \eqref{eq:formel-rho-tilde} follows from
\eqref{eq:def_tilde_eta} and  \(\widetilde\eta''(0)=0\).
This finishes the proof of Theorem~\ref{thm:eq-rho-tilde}.
\qed

\section{Regularity of \(h\) and \(\widetilde h\)}
\label{app:reg}

In this section we prove the statements in
Theorem~\ref{thm:eq-rho-tilde} on the regularity of the functions
\(h\) and \(\widetilde h\).
More precisely, we prove the following:
\begin{prop}
 \label{prop:reg}
 Let \(\psi\) satisfy \eqref{eq:eigen} and let
 \(h\) be as in \eqref{eq:rho_first}:
\begin{align}
   \label{eq:h}
   h(x)&=
  \int_{\R^{3(N-1)}}
  |\nabla\psi|^{2} 
  \,dx_{2}\cdots dx_{N}
  \nonumber\\
  &\quad-\sum_{j=2}^{N}\int_{\R^{3(N-1)}}\frac{Z}{|x_{j}|}\,\psi^{2}
  \,dx_{2}\cdots dx_{N}
  \nonumber\\ 
  &\quad+\sum_{1\leq j<k\leq N}\int_{\R^{3(N-1)}}
  \frac{1}{|x_{j}-x_{k}|}\,\psi^{2}
  \,dx_{2}\cdots dx_{N},
\end{align}
and \(\widetilde h\) as in \eqref{def_h_tilde}: 
\begin{align}
  \widetilde h(r)&=\int_{\mathbb{S}^{2}}h(r\omega)\,d\omega.
  \nonumber
\end{align}
Then \(h\in C^{\alpha}_{\text{{\rm loc}}}(\R^{3}\setminus\{0\})\) and 
\(\widetilde h\in C^{0}([0,\infty))\cap 
C^{\alpha}_{\text{{\rm loc}}}((0,\infty))\) 
for all \(\alpha\in(0,1)\). 
\end{prop}
\begin{remark}
  From the proof of Proposition~\ref{prop:reg} follows 
  the continuity of the function \(t_{1}\) in 
  Remark~\ref{t_1}.
\end{remark}
\begin{proof}
For convenience, we shall often write \(\int\equiv\int_{\R^{3(N-1)}}\).
Let
\begin{align}
  \label{eq:def_j}
  J_{1}(x)&=\int|\nabla\psi|^{2} 
  \,dx_{2}\cdots dx_{N},
  \nonumber\\
  J_{2}(x)&=\sum_{j=2}^{N}\int\frac{Z}{|x_{j}|}\psi^{2}
  \,dx_{2}\cdots dx_{N},
  \nonumber\\
  J_{3}(x)&=\sum_{1\leq j<k\leq N}\int\frac{1}{|x_{j}-x_{k}|}\psi^{2}
  \,dx_{2}\cdots dx_{N}.
\end{align}
For the proof of the regularity of the functions \(J_{1}, J_{2}, J_{3}\) we shall make use of the following lemmas.
The proof of the first lemma is trivial, using 
\(|(x_{1},x_{2},\ldots,x_{N})|\geq|(x_{2},\ldots,x_{N})|\). 
\begin{lemma}
  \label{lem:one}
  Let \(\alpha\in(0,1)\) and \(x_{0}\in\R^{3}\).
  Assume that the real function \(G=G(x_{1},\ldots,x_{N})\) satisfies: For
  all \(R>0\) there exists constants \(C,\gamma\) such that
  \begin{align}
    \label{eq:cond1}
    &\sup_{x,y\in B(x_{0},R)}
     \frac{|G(x,x_{2},\ldots,x_{N})-G(y,x_{2},\ldots,x_{N})|}{|x-y|^{\alpha}} 
     \leq
     \nonumber\\ 
     &C\,\exp\big(-\gamma|(x_{0},x_{2},\ldots,x_{N})|\big)
     \quad\text{for all}\quad(x_{0},x_{2},\ldots,x_{N})\in\R^{3N}.
  \end{align}
  Then the function
  \begin{align}
    \eta(x)\equiv
    \int_{\R^{3(N-1)}}
    \!\!\!\!\!\!\!\!\!\!
    G(x,x_{2},\ldots,x_{N})\,dx_{2}\cdots dx_{N}
    \nonumber
  \end{align}
  is in \(C^{\alpha}_{\text{loc}}(\R^{3})\).
\end{lemma}

We next prove the following lemma.
\begin{lemma}
  \label{lem:two}
  Let \(\alpha\in(0,1)\).
  Assume that the real valued function \(K=K(x_{1},\ldots,x_{N})\) satisfies  
  \eqref{eq:cond1} and that there exists constants \(C,\gamma\) such that
  \begin{align}
    \label{eq:cond2}
    &|K(x_{1},\dots,x_{N})|
    \nonumber\\
    &\leq C\,\exp\big(-\gamma|(x_{1},\ldots,x_{N})|\big)
    \text{ for all } (x_{1},\ldots,x_{N})\in\R^{3N}.
  \end{align}
  Then: 

(a)  

  For all \(j,k\in\{1,\ldots,N\},j\neq k\), the function 
  \begin{align}
    \zeta(x_{1})\equiv\int_{\R^{3(N-1)}}\frac{1}{|x_{j}-x_{k}|}\,
    K(x_{1},\ldots,x_{N})\,dx_{2}\cdots dx_{N}
    \nonumber
  \end{align}
    is in \(C^{\alpha}_{\text{{\rm loc}}}(\R^{3})\).

(b)  

  For \(j\geq2\) the function 
   \begin{align}
    \mu(x_{1})\equiv\int_{\R^{3(N-1)}}\frac{1}{|x_{j}|}\,
    K(x_{1},\ldots,x_{N})\,dx_{2}\cdots dx_{N}
    \nonumber
  \end{align}
  is in \(C^{\alpha}_{\text{{\rm loc}}}(\R^{3})\).
\end{lemma}
\begin{proof}
  Assume first that \(j\neq1\neq k\). Let \(x,y\in B(x_{0},R)\), 
  then by \eqref{eq:cond1},
  \begin{align}
    &\frac{|\zeta(x)-\zeta(y)|}{|x-y|^{\alpha}}
    \leq 
    \nonumber\\
    &\int\frac{1}{|x_{j}-x_{k}|}
    \frac{|K(x,x_{2},\ldots,x_{N})-K(y,x_{2},\ldots,x_{N})|}{|x-y|^{\alpha}} 
    \,dx_{2}\cdots dx_{N}
    \nonumber\\
    &\leq C \int\frac{1}{|x_{j}-x_{k}|}
    \exp\big(-\gamma|(x_{0},x_{2},\ldots,x_{N})|\big)\,dx_{2}\cdots
    dx_{N}.
    \nonumber
  \end{align}
  By equivalence of norms in \(\R^{3N}\) there is a constant \(c_{0}\)
  such that
  \begin{align}
    &\frac{|\zeta(x)-\zeta(y)|}{|x-y|^{\alpha}}\leq
    C\Big(\prod_{l=2,l\neq j,k}^{N}
    \int_{\R^{3}}\exp(-\gamma c_{0}|x_{l}|)\,dx_{l}\Big)\times
    \nonumber\\
    &\quad
    \times\int_{\R^{6}}
    \frac{1}{|x_{j}-x_{k}|}\exp(-\gamma c_{0}(|x_{j}|+|x_{k}|))
    \,dx_{j}\,dx_{k}\leq C,\,\,x,y\in B(x_{0},R).
    \nonumber
  \end{align}
 The last inequality is an application of the following inequality (with
 \(n=3, \lambda=1, p=r=6/5\)):
 (see Lieb and Loss~\cite[Theorem 4.3]{Li-Loss})
 
\vskip 0.4cm
\noindent
{\bf Hardy-Littlewood-Sobolev Inequality:}
{\emph  {Let \(p,r>1\) and \(0<\lambda<n\) with \(1/p+\lambda/n+1/r=2\). Let
  \(f\in L^{p}(\R^{n})\) and \(h\in L^{r}(\R^{n})\). Then there exists 
  a sharp constant \(C(n,\lambda,p)\), independent of \(f\) and \(h\), 
  such that
  \begin{align}
    \nonumber
    \left|\int_{\R^{n}}\int_{\R^{n}}
     f(x)|x-y|^{-\lambda}h(y)\,dx\,dy\right|
     \leq C(n,\lambda,p)\|f\|_{p}\|h\|_{r}.
    \nonumber
  \end{align}
 }}
  This  proves Lemma~\ref{lem:two} (a) when  \(j\neq1\neq k\). 

  Assume now that \(j=1\). We assume without loss that \(k=2\).
  Then, with \(x_{1},\bar x_{1}\in B(x_{0},R)\),
  \begin{align}
    &\frac{|\zeta(x_{1})-\zeta(\bar x_{1})|}{|x_{1}-\bar x_{1}|^{\alpha}}\leq
    \nonumber\\
    &{}\int\frac{1}{|x_{1}-x_{2}|}
    \frac{|K(x_{1},x_{2},\ldots,x_{N})
   -K(\bar x_{1},x_{2},\ldots,x_{N})|}{|x_{1} -\bar x_{1}|^{\alpha}} 
    \,dx_{2}\cdots dx_{N}
    \nonumber\\
    &+\int \left|\frac{1}{|x_{1}-x_{2}|}-\frac{1}{|\bar x_{1}-x_{2}|}\right|
    |x_{1}-\bar x_{1}|^{-\alpha}|K(\bar x_{1},x_{2},\ldots,x_{N})|
    \,dx_{2}\cdots dx_{N}
    \nonumber\\
    &\equiv (\text{A})+(\text{B}).
    \nonumber
  \end{align}
  For \((\text{A})\), using~\eqref{eq:cond1} and equivalence of 
  norms in \(\R^{3N}\), we get
  \begin{align}
    (\text{A})&\leq
    C\int\frac{1}{|x_{1}-x_{2}|}\exp\big(-\gamma c_{0}(|x_{0}|+|x_{2}|
    +\cdots+|x_{N}|)\big)\,dx_{2}\cdots dx_{N}
    \nonumber\\
    &\leq C\int_{\R^{3}}\frac{1}{|x_{1}-x_{2}|}
    \exp(-\gamma c_{0}|x_{2}|)\,dx_{2}\leq C(x_{0},R).
    \nonumber
  \end{align}
  As for \((\text{B})\),
  we apply the following inequality; for the convenience of the reader, we
  give the proof (borrowed from Lieb and Loss~\cite[(3) p. 225]{Li-Loss}).

\vskip 0.2cm

For \(\alpha\in(0,1)\): 
 \begin{align}
    \label{lem:Li-Loss}
    &\left| \frac{1}{|x-z|}-\frac{1}{|y-z|}\right| 
    |x-y|^{-\alpha}
    \leq  
     \nonumber\\
    &|x-z|^{-1-\alpha}+|y-z|^{-1-\alpha} \quad\text{for all}\quad
     x,y,z\in\R^{3}.
  \end{align}

  {\emph {Proof}} of \eqref{lem:Li-Loss}.
    By H\"older's inequality we have, for \(b>1\), \(\alpha\in(0,1)\), 
    \begin{align}
      \nonumber
       1-b^{-1}=\int_{1}^{b}t^{-2}\,dt
       \leq\left(\int_{1}^{b} \,dt\right)^{\alpha}
       \left(\int_{1}^{\infty}t^{-2/(1-\alpha)}\,dt\right)^{1-\alpha}
       \leq(b-1)^{\alpha}.
    \end{align}
    Substituting \(b/a\) for \(b\), with \(a>0\), this gives
    \begin{align}
      \nonumber
      \left|b^{-1}-a^{-1}\right|
      \leq |b-a|^{\alpha}\max\{a^{-1-\alpha},b^{-1-\alpha}\}.
    \end{align}
    So, for \(x,y,z\in\R^{3}\), using \(\max\{s,t\}\leq s+t\) and the 
    triangle inequality in \(\R^{3}\), we have 
    \begin{align}
      \left| |x-z|^{-1}-|y-z|^{-1}\right|
      \leq |x-y|^{\alpha}\left\{
      |x-z|^{-1-\alpha}+|y-z|^{-1-\alpha}\right\}.
      \nonumber
    \end{align}
  \qed

  In this way, by \eqref{eq:cond2} and equivalence of norms in
  \(\R^{3N}\),
  \begin{align}
    (\text{B})&\leq
    \int
    \left(\frac{|K(\bar x_{1},x_{2},\ldots,x_{N})|}{|x_{1}-x_{2}|^{1+\alpha}}
    +\frac{|K(\bar x_{1},x_{2},\ldots,x_{N})|}{|\bar
      x_{1}-x_{2}|^{1+\alpha}}  
    \right)
    \,dx_{2}\cdots dx_{N}
    \nonumber\\ 
    &\leq
    C\prod_{j=3}^{N}
    \int_{\R^{3}}\exp(-\gamma c_{0}|x_{j}|)\,dx_{j}
    \nonumber\\
    &\qquad\times
    \int_{\R^{3}}\left(\frac{1}{|\bar x_{1}-x_{2}|^{1+\alpha}}
    +\frac{1}{|x_{1}-x_{2}|^{1+\alpha}}\right)\exp(-\gamma
    c_{0}|x_{2}|)\,dx_{2} 
    \nonumber\\
    &\leq C(x_{0},R),
    \nonumber
  \end{align}
  since \(x_{1},\bar x_{1}\in B(x_{0},R)\). This finishes the proof of 
  Lemma~\ref{lem:two} (a).

  The proof of \((b)\) is similar to that of \((a)\) so we omit the details.
\end{proof}
The proof of the following fact is straightforward:

There exist constants \(C=C(\gamma,R)\) and
  \(\widetilde\gamma=\widetilde\gamma(\gamma)\) such that 
  \begin{align}
    \label{lem:three}
    \exp(-\gamma|(x,\ldots,x_{N})|)
    \leq C\,\exp(-\widetilde\gamma|(x_{0},\ldots,x_{N})|)
  \end{align}
for all \(x\in B(x_{0},R)\).

Using  this and Lemma~\ref{lem:one}, \ref{lem:two}, 
we shall prove the following lemma on the regularity of the functions 
\(J_{1}, J_{2}\) and \(J_{3}\) from \eqref{eq:def_j}.
\begin{lemma} 
  \label{lem:reg_Js}
  Let \(J_{1}, J_{2}\) and \(J_{3}\) be as in  \eqref{eq:def_j}.
  Then 

  \begin{enumerate}
  \item 
    \(J_{2}, J_{3}\in C_{\text{{\rm loc}}}^{\alpha}(\R^{3})\)
    for all \(\alpha\in(0,1)\).
  \item
     \(J_{1}\in C_{\text{{\rm
          loc}}}^{\alpha}(\R^{3}\setminus\{0\})\)
         for all \(\alpha\in(0,1)\).
  \end{enumerate}
\end{lemma}

Herefrom follow the regularity properties of the function \(h\) stated in
Proposition~\ref{prop:reg}. 

{\emph{Proof}} of Lemma~\ref{lem:reg_Js} (i).

Firstly, by Theorem~\ref{thm:grad_est} and Remark~\ref{rem:1}, 
\begin{align}
  \label{eq:both_exp}
  |\psi({\bf x})|
  \,\,,\,\,\,
  |\nabla\psi({\bf x})|\leq C\exp(-\gamma|{\bf x}|)\,\,,
  \,\,\,{\bf x}\in\R^{3N},
\end{align}
which gives \eqref{eq:cond2} for \(K=\psi^{2}\).

Next we verify that for \(G=\psi^{2}\), \eqref{eq:cond1} is
fulfilled. Then Lemma~\ref{lem:two} can be applied with
\(K=\psi^{2}\), and the H\"older-continuity of \(J_{2}\) and \(J_{3}\) 
follows.

Given \(x_{0}\in\R^{3}, R>0, \alpha\in(0,1)\), and \(x,y\in B(x_{0},R)\).
Using that (see e.\ g.\ Mal\'y and Ziemer~\cite[Theorem 1.41]{Maly-Ziemer}) 
(here, \((x,x_{2},\ldots,x_{N})=(x,x')\), \(x'\in\R^{3(N-1)}\))
\begin{align}
  \psi^{2}(x,x')
  &-\psi^{2}(y,x')
  =\int_{0}^{1}\frac{\partial}{\partial s}
  \big[\psi^{2}(sx+(1-s)y,x')\big]\,ds
  \label{eq:diff_grad}
  \nonumber\\ 
  &=\int_{0}^{1}\big[\nabla_{1}(\psi^{2})
  (sx+(1-s)y,x')\big]\cdot(x-y)\,ds 
\end{align}
and that \(sx+(1-s)y\in B(x_{0},R)\) for all
\(s\in[0,1]\) we get, with \eqref{eq:both_exp} and \eqref{lem:three},
\begin{align}
  &\frac{|\psi^{2}(x,x')-\psi^{2}(y,x')|}{|x-y|^{\alpha}}
  \nonumber\\
  &\leq 2
  |x-y|^{1-\alpha}\int_{0}^{1}|\nabla_{1}\psi(sx+(1-s)y,x')|
  \cdot|\psi(sx+(1-s)y,x')|\,ds
  \nonumber\\
  &\leq 2
  (2R)^{1-\alpha}\int_{0}^{1}C\,
  \exp(-2\gamma|(sx+(1-s)y,x')|)\,ds
  \nonumber\\
  &\leq C\,\exp(-\widetilde\gamma|(x_{0},\dots,x_{N})|),
  \label{eq:midl}
 \end{align}
so \eqref{eq:cond1} follows for \(\alpha\in(0,1)\). 

This proves (i) of Lemma~\ref{lem:reg_Js}.

To prove (ii), we write \(\psi\) as in the proof of
Theorem~\ref{thm:grad_est}: \(\psi=e^{F-F_{1}}\psi_{1}\), with \(F\)
  and \(F_{1}\) as in \eqref{eq:def_F} and \eqref{def_F_1}. Then 
  \begin{align}
    \label{eq:int_grad_psi}
    J_{1}&(x)=\int|\nabla\psi|^{2}\,dx'
    =\int|\nabla F|^{2}\psi^{2}\,dx'
    +\int|\nabla F_{1}|^{2}\psi^{2}\,dx'
    \nonumber\\
    &-2\int\big(\nabla F\cdot\nabla F_{1}\big)\psi^{2}\,dx'
    +\int e^{2(F-F_{1})}|\nabla\psi_{1}|^{2}\,dx'
    \nonumber\\
    &+2\int\big(\nabla F\cdot\nabla\psi_{1}\big)
    e^{2(F-F_{1})}\psi_{1}\,dx'
    -2\int\big(\nabla F_{1}\cdot\nabla\psi_{1}\big)
    e^{2(F-F_{1})}\psi_{1}\,dx'
    \nonumber\\
    &\equiv
    I_{1}(x)+I_{2}(x)+I_{3}(x)+I_{4}(x)+I_{5}(x)+I_{6}(x).
  \end{align}
Using the idea from \eqref{eq:diff_grad} twice (on \(|\nabla
F_{1}|^{2}\) and \(\psi^{2}\), respectively), the estimates 
\eqref{eq:est_F}, \eqref{eq:both_exp}, \eqref{eq:midl},
and \eqref{lem:three}, we have, with
\(x_{0}\in\R^{3}, R>0, \alpha\in(0,1)\), and \(x,y\in B(x_{0},R)\):
\begin{align}
  &\frac{\big||\nabla F_{1}|^{2}\psi^{2}(x,x')-
   |\nabla F_{1}|^{2}\psi^{2}(y,x')\big|}{|x-y|^{\alpha}}
  \leq
  \nonumber\\
  &\qquad\quad
  \frac{\big||\nabla F_{1}(x,x')|^{2}-|\nabla
  F_{1}(y,x')|^{2}\big|}{|x-y|^{\alpha}}|\psi(x,x')|^{2}
  \nonumber\\
  &\qquad\quad
  +|\nabla F_{1}(y,x')|^{2}\frac{\big|\psi^{2}(x,x')
  -\psi^{2}(y,x'x)\big|}{|x-y|^{\alpha}}
  \nonumber\\
  &\quad\quad\leq C\,|x-y|^{1-\alpha}\big\|\nabla\big(|\nabla
  F_{1}|^{2}\big)\big\|_{\infty}\exp(-2\gamma|(x,x_{2},\ldots,x_{N})|)
  \nonumber\\
  &\qquad\quad\quad
  +2\widetilde C|x-y|^{1-\alpha}\big\||\nabla F_{1}|^{2}\big\|_{\infty}
  \exp(-\widetilde\gamma|(x_{0},\ldots,x_{N})|)
  \nonumber\\
  &\quad\quad
  \leq \bar C\,\exp(-\bar\gamma(x_{0},\ldots,x_{N})|).
  \nonumber
\end{align}
By Lemma~\ref{lem:one}, with \(G=|\nabla
F_{1}|^{2}\psi^{2}(x_{1},\cdots,x_{N})\),  
this implies that \(I_{2}\in C_{\text{loc}}^{\alpha}(\R^{3})\). 

Using the same ingredients, writing
\(\nabla\psi_{1}=\nabla(e^{F_{1}-F}\psi)\), gives \eqref{eq:cond1}
and \eqref{eq:cond2} with 
\begin{align}
  G&=K=e^{2(F-F_{1})}|\nabla\psi_{1}|^{2}\nonumber\\
\intertext{and} 
  G&=K=\big(\nabla F_{1}\cdot
  \nabla\psi_{1}\big)e^{2(F-F_{1})}\psi_{1},
  \nonumber
\end{align}
and so by Lemma~\ref{lem:one}, 
\(I_{4},I_{6}\in C_{\text{loc}}^{\alpha}(\R^{3})\).

The remaining terms are those involving the function \(F\),
namely \(I_{1}\), \(I_{3}\) and \(I_{5}\).

Note that 
\begin{align}
  \label{eq:grad_F}
  &\nabla F(x_{1},\ldots,x_{N})
  ={}-\frac{Z}{2}\Big(\frac{x_{1}}{|x_{1}|},\ldots,\frac{x_{N}}{|x_{N}|}\Big)
  \nonumber\\
  &+\frac{1}{4}\left(
   \sum_{k=2}^{N}\frac{x_{1}-x_{k}}{|x_{1}-x_{k}|}
   \,,\ldots,\,
   \sum_{k=1,k\neq j}^{N}\frac{x_{j}-x_{k}}{|x_{j}-x_{k}|}
   \,,\ldots,\,
   \sum_{k=1}^{N-1}\frac{x_{N}-x_{k}}{|x_{N}-x_{k}|}
   \right)
\end{align}
and so
\begin{align}
  |\nabla F(x_{1},\ldots,x_{N})|^{2}
  &=\frac{NZ^{2}}{4}-\frac{Z}{8}
  \sum_{j,k=1,k\neq j}^{N}\frac{x_{j}}{|x_{j}|}\cdot
  \frac{x_{j}-x_{k}}{|x_{j}-x_{k}|}
  \nonumber\\
  &\quad+\frac{1}{16}\sum_{j,k,l=1,k\neq j,l\neq j}^{N}
  \frac{x_{j}-x_{k}}{|x_{j}-x_{k}|}\cdot\frac{x_{j}-x_{l}}{|x_{j}-x_{l}|}.
  \nonumber
\end{align}
In this way,
\begin{align}
  &I_{1}(x_{1})=
  \frac{NZ^{2}}{4}\int\psi^{2}
  \,dx_{2}\cdots dx_{N}
  \nonumber\\
  &\qquad\qquad
  -\frac{Z}{8}\sum_{j,k=1,k\neq j}^{N}
  \int\frac{x_{j}}{|x_{j}|}\cdot\frac{x_{j}-x_{k}}{|x_{j}-x_{k}|}\,\psi^{2}
  \,dx_{2}\cdots dx_{N}
  \nonumber\\
  &\qquad\qquad
  +\frac{1}{16}\sum_{j,k,l=1,k\neq j,l\neq j}^{N}
  \frac{x_{j}-x_{k}}{|x_{j}-x_{k}|}\cdot\frac{x_{j}-x_{l}}{|x_{j}-x_{l}|}
  \,\psi^{2}
  \,dx_{2}\cdots dx_{N}
  \nonumber\\
  &
  =\frac{NZ^{2}}{4}
  \rho(x_{1})-\frac{Z}{8}\sum_{j,k=1,k\neq j}^{N}\kappa_{j,k}(x_{1})
  +\frac{1}{16}\sum_{j,k,l=1,k\neq j,l\neq j}^{N}
  \nu_{j,k,l}(x_{1}).
  \label{eq:I_1}
\end{align}
Note that \(\nu_{j,k,l}=\nu_{j,l,k}\). 

Using the ideas above, Lemma~\ref{lem:one} implies 
that the following functions from \eqref{eq:I_1} (with the
mentioned choices of \(G\) satisfying \eqref{eq:cond1}) 
are all in \(C_{\text{loc}}^{\alpha}(\R^{3})\):
\begin{align}
  \rho&: \quad G=\psi^{2},
  \nonumber\\  
  \kappa_{j,k} , \quad j,k\neq1, j\neq k&:\quad
  G=\frac{x_{j}}{|x_{j}|}\cdot\frac{x_{j}-x_{k}}{|x_{j}-x_{k}|}\,\psi^{2},
  \nonumber\\
  \nu_{j,k,k},j\neq k&:\quad
  G=\frac{x_{j}-x_{k}}{|x_{j}-x_{k}|}
  \cdot\frac{x_{j}-x_{k}}{|x_{j}-x_{k}|}\,\psi^{2}=\psi^{2},
  \nonumber\\
  \nu_{j,k,l}, \quad j,k,l\neq1,l\neq j\neq k&:\quad
  G=\frac{x_{j}-x_{k}}{|x_{j}-x_{k}|}  
  \cdot\frac{x_{j}-x_{l}}{|x_{j}-x_{l}|}\,\psi^{2}.
  \nonumber
\end{align}
Likewise, Lemma~\ref{lem:two} implies (with the mentioned choices of
\(G=K\) satisfying \eqref{eq:cond1} and \eqref{eq:cond2}) 
that the following functions from \eqref{eq:I_1} are all in
\(C_{\text{loc}}^{\alpha}(\R^{3})\): 
\begin{align}
  \kappa_{j,1}, \quad j\neq1&:\quad
  G=K=\frac{x_{j}\cdot(x_{j}-x_{1})}{|x_{j}|}\,\psi^{2},
  \nonumber\\
  \nu_{j,1,l},\quad j,l\neq1, j\neq l&:\quad
  G=K=\frac{(x_{j}-x_{1})\cdot(x_{j}-x_{l})}{|x_{j}-x_{l}|}\,\psi^{2},
  \nonumber
\end{align}
From the decomposition of \(I_{1}\) in \eqref{eq:I_1} we are left with
\begin{align}
  \label{eq:not_alpha}
  \kappa_{1,k}(x_{1})=
  &\int\frac{x_{1}}{|x_{1}|}\cdot\frac{x_{1}-x_{k}}{|x_{1}-x_{k}|}\,\psi^{2}
  \,dx'
  \quad,\quad k=2,\ldots,N,
  \intertext{ and } 
  \label{eq:two_terms}
  \nu_{1,k,l}(x_{1})=
  &\int\frac{x_{1}-x_{k}}{|x_{1}-x_{k}|}\cdot
  \frac{x_{1}-x_{l}}{|x_{1}-x_{l}|}\,\psi^{2}
  \,dx'
  \quad,\quad
  k,l\in\{2,\ldots,N\}, k\neq l.
 \end{align}
 Note that
 \begin{align}
   &\int\frac{x_{1}}{|x_{1}|}\cdot\frac{x_{1}-x_{k}}{|x_{1}-x_{k}|}\,\psi^{2}
   \,dx_{2}\cdots dx_{N}   
   \nonumber\\
   &\qquad=\frac{1}{|x_{1}|}\int
   \frac{1}{|x_{1}-x_{k}|}
   \big(x_{1}\cdot(x_{1}-x_{k})\,\psi^{2}\big)\,dx_{2}\cdots dx_{N}.
   \nonumber
 \end{align}
 The function \(1/|x_{1}|\) is smooth for \(x_{1}\neq0\)
 and therefore in
 \(C^{\alpha}_{\text{{\rm loc}}}(\R^{3}\setminus\{0\})\). The function
  \(x_{1}\cdot(x_{1}-x_{k})\,\psi^{2}\) satisfies
\eqref{eq:cond1} 
and \eqref{eq:cond2} (by the same ideas as above), so Lemma~\ref{lem:two} (a)
implies that the function
\begin{align}
  \nonumber
  \int\frac{1}{|x_{1}-x_{k}|}
  \big(x_{1}\cdot(x_{1}-x_{k})\,\psi^{2}\big)\,dx_{2}\cdots dx_{N}
\end{align}
is in \(C^{\alpha}_{\text{{\rm loc}}}(\R^{3})\).  The functions 
 in \eqref{eq:not_alpha} are therefore in
\(C^{\alpha}_{\text{{\rm loc}}}(\R^{3}\setminus\{0\})\).  

As for the functions in \eqref{eq:two_terms}, these are all in
\(C^{\alpha}_{\text{{loc}}}(\R^{3})\), which can be seen by applying
the previous ideas, in particular Lemma~\ref{lem:Li-Loss}, \eqref{lem:three},
\eqref{eq:both_exp} and \eqref{eq:midl}.

This proves that 
\(I_{1}\in C^{\alpha}_{\text{{\rm loc}}}(\R^{3}\setminus\{0\})\).

As for \(I_{3}\) (see \eqref{eq:int_grad_psi} and \eqref{eq:grad_F}),
with \(\nabla=(\nabla_{1},\ldots,\nabla_{N})\), 
\begin{align}
  \label{eq:I3}
  I_{3}(x)&=
  Z\sum_{j=1}^{N}\int\Big(\frac{x_{j}}{|x_{j}|}\cdot\nabla_{j}F_{1}\Big)
  \psi^{2}\,dx_{2}\cdots dx_{N}
   \nonumber\\
   &\quad-\frac{1}{2}\sum_{j,k=1,j\neq k}^{N}\int
  \Big(\frac{x_{j}-x_{k}}{|x_{j}-x_{k}|}\cdot\nabla_{j}F_{1}\Big)
  \psi^{2}\,dx_{2}\cdots dx_{N}.
\end{align}
The terms in the first sum with \(j\neq1\) are in 
\(C_{\text{loc}}^{\alpha}(\R^{3})\), due to 
Lemma~\ref{lem:two} (b), with
\(G=K=(x_{j}\cdot\nabla_{j}F_{1})\,\psi^{2}\)
satisfying \eqref{eq:cond1} and \eqref{eq:cond2}. (To see this, use 
the previous ideas; to apply the idea from
\eqref{eq:diff_grad} to \(\nabla_{j}F_{1}\) we use that \(F_{1}\) is smooth).
The terms in the second sum in \eqref{eq:I3} are all in
\(C_{\text{loc}}^{\alpha}(\R^{3})\), due to Lemma~\ref{lem:two} (a),
applied with
\(G=K=\big((x_{j}-x_{k})\cdot
\nabla_{j}F_{1}\big)\,\psi^{2}\).
The term with \(j=1\) is in
\(C_{\text{loc}}^{\alpha}(\R^{3}\setminus\{0\})\). This can be seen
by following the ideas in the proof of the regularity properties of the
function in \eqref{eq:not_alpha}, now
using Lemma~\ref{lem:one} with \(G=(x_{1}\cdot\nabla_{1}F_{1})\psi^{2}\).

The statements and proofs are similar for 
\begin{align}
  \label{eq:I_5}
  I_{5}(x)&=
  -Z\sum_{j=1}^{N}\int\Big(\frac{x_{j}}{|x_{j}|}\cdot\nabla_{j}\psi_{1}\Big)
  e^{2(F-F_{1})}\psi_{1}\,dx_{2}\cdots dx_{N}
  \nonumber\\
  &{}\quad
  +\frac{1}{2}\sum_{j,k=1,j\neq k}^{N}\int
  \Big(\frac{x_{j}-x_{k}}{|x_{j}-x_{k}|}\cdot\nabla_{j}\psi_{1}\Big)
  e^{2(F-F_{1})}\psi_{1}\,dx_{2}\cdots dx_{N}.
\end{align}
That is, the functions in the first sum in \eqref{eq:I_5} with
\(j\geq2\) and those in the second sum are all in 
\(C^{\alpha}_{\text{{\rm loc}}}(\R^{3})\), whereas the function in
the first sum with \(j=1\) is only in 
\(C^{\alpha}_{\text{{\rm loc}}}(\R^{3}\setminus\{0\})\). 
To prove this we use the inequality
(with \({\bf x}=(x_{1},\ldots,x_{N})\)):
 \begin{align}
    \nonumber
    |\psi_{1}|_{C^{1,\alpha}(B({\bf x},R/2))}
    \leq C\, \sup_{y\in(B({\bf x},R))}|\psi_{1}(y)|
    \leq C\,\exp(-\gamma|(x_{1},\ldots,x_{N})|). 
  \end{align}
This inequality follows from \eqref{eq:est_F}, \eqref{eq:apriori_1} and
\eqref{eq:both_exp} (remember that \(\psi_{1}=e^{F_{1}-F}\psi\)).

This proves that \(I_{5}\in
C^{\alpha}_{\text{{loc}}}(\R^{3}\setminus\{0\})\), and so finishes the 
proof that \(J_{1}\in C^{\alpha}_{\text{{\rm loc}}}(\R^{3}\setminus\{0\})\).
(See \eqref{eq:int_grad_psi}). This proves (ii) and therefore
Lemma~\ref{lem:reg_Js}.
\qed

  That \(\widetilde h\in C^{\alpha}_{\text{{\rm loc}}}((0,\infty))\)
  is a consequence of the foregoing and of the following proposition:
  \begin{prop}
    Assume \(f\in C^{\alpha}_{\text{{\rm
          loc}}}(\R^{3}\setminus\{0\})\), \(\alpha\in(0,1)\).
    Then \(\widetilde f\in
    C^{\alpha}_{\text{{\rm loc}}}((0,\infty))\), where
    \begin{align}
      \nonumber
      \widetilde f(r)=\int_{\sphere^{2}}f(r\omega)\,d\omega.
    \end{align}
  \end{prop}
  \begin{proof}
    Let \(r\in(0,\infty)\). For all 
  \(x_{0}\in A=\{x\in\R^{3}\,|\,|x|=r\}\), 
    choose \(R=R(x_{0})\) and 
    \(C=C(x_{0})\) such that
      \begin{align}
        \label{eq:Holder_f}
        \sup_{x,y\in B(x_{0},R(x_{0}))} 
        \frac{|f(x)-f(y)|}{|x-y|^{\alpha}}\leq C(x_{0}).
      \end{align}
    This is possible, since 
    \(f\in C^{\alpha}_{\text{{\rm loc}}}(\R^{3}\setminus\{0\})\).
     Then 
    \begin{align}
      \nonumber
      A\subset\bigcup_{x_{0}\in A}B(x_{0},R(x_{0})). 
    \end{align}
    Using compactness of \(A\), choose \(x_{1},\ldots,x_{m}\in A\)
    such that 
    \begin{align}
      \nonumber
      A\subset\bigcup_{j=1}^{m}B(x_{j},R(x_{j})). 
    \end{align}
    Choose \(\epsilon\in(0,r)\) such that
    \begin{align}
      \nonumber 
     \{y\in\R^{3}\,|\,r-\epsilon<|y|<r+\epsilon\}\subset
      \bigcup_{j=1}^{m}B(x_{j},R(x_{j})). 
    \end{align}
    Then, for all \(s,t\in(r-\epsilon,r+\epsilon)\) and all
    \(\omega\in\sphere^{2}\) there exists \(j\in\{1,\ldots,m\}\) such
    that \(s\omega,t\omega\in B(x_{j},R(x_{j}))\) and therefore by
    \eqref{eq:Holder_f}, 
    \begin{align}
      \frac{|f(s\omega)-f(t\omega)|}{|s-t|^{\alpha}}
      =\frac{|f(s\omega)-f(t\omega)|}{|s\omega-t\omega|^{\alpha}}
      \leq C(x_{j}).
      \nonumber
    \end{align}
    So with \(C=\max\{C(x_{1}),\ldots,C(x_{m})\}\),
    \begin{align}
      \frac{|f(s\omega)-f(t\omega)|}{|s-t|^{\alpha}}
      \leq C,
      \text{ for all } s,t\in (r-\epsilon,r+\epsilon)\,\, \text{and all }
      \omega\in\sphere^{2}.
      \nonumber
    \end{align}

    This implies that
    \begin{align}
      &\frac{|\widetilde f(s)-\widetilde f(t)|}{|s-t|^{\alpha}}
      =\frac{|\int_{\sphere^{2}}
             (f(s\omega)-f(t\omega))\,d\omega|
            }{|s-t|^{\alpha}}
      \nonumber\\      
      &\quad\leq \int_{\sphere^{2}}
      \frac{|f(s\omega)-f(t\omega)|}{|s\omega-t\omega|^{\alpha}}\,d\omega 
      \leq C,\quad
      \text{ for all } s,t\in(r-\epsilon,r+\epsilon).
      \nonumber
    \end{align}
     This proves that \(\widetilde f\in
    C^{\alpha}_{\text{{\rm loc}}}((0,\infty))\).
  \end{proof}

To prove that \(\widetilde h\in C^{0}([0,\infty))\), we apply the
following:
\begin{prop}
  \label{prop:cont_int}
  Assume \(f\in C^{\alpha}_{\text{{\rm loc}}}(\R^{3})\). Then \(\widetilde 
  f\in C^{0}([0,\infty))\), where
  \begin{align}
    \nonumber
    \widetilde f(r)=\int_{\sphere^{2}}f(r\omega)\,d\omega.
  \end{align}
\end{prop}
\begin{proof}
  The function \(f\) is continuous in \(\R^{3}\), since it is in 
  \(C^{\alpha}_{\text{{\rm loc}}}(\R^{3})\). Let
  \(r\in[0,\infty)\). Then
  \begin{align}
    \lim_{s\to r}f(s\omega)=f(r\omega)\quad 
    \text{ for all } \omega\in\sphere^{2}.
    \nonumber
  \end{align}
  Using the supremum of \(f\) on a sufficiently large compact set in
  \(\R^{3}\) as a dominant, Lebesque's Dominated Convergence Theorem
  gives us that
  \begin{align}
    \lim_{s\to r}\int_{\sphere^{2}}f(s\omega)\,d\omega    
    =\int_{\sphere^{2}}f(r\omega)\,d\omega.
    \nonumber
  \end{align}
  Therefore \(f\in C^{0}([0,\infty))\).
\end{proof}
Recall the proof of the fact that 
\(h\in C^{\alpha}_{\text{{\rm loc}}}(\R^{3}\setminus\{0\})\).
In fact, the only terms in the decomposition of \(h\)  (see
\eqref{eq:h}, \eqref{eq:def_j}, \eqref{eq:int_grad_psi}, and 
\eqref{eq:grad_F}) that are only in
\(C^{\alpha}_{\text{{\rm loc}}}(\R^{3}\setminus\{0\})\) and not in
\(C^{\alpha}_{\text{{\rm loc}}}(\R^{3})\) are the functions
\begin{align}
  \label{eq:except}
  &\int\frac{x_{1}}{|x_{1}|}\cdot\frac{x_{1}-x_{k}}{|x_{1}-x_{k}|}\,\psi^{2}
  \,dx_{2}\cdots dx_{N}
  \quad,\quad k=2,\ldots,N,
  \nonumber\\
  &\int\Big(\frac{x_{1}}{|x_{1}|}\cdot\nabla_{1}F_{1}\Big)\psi^{2}
  \,dx_{2}\cdots dx_{N},
  \nonumber\\
  &\int\Big(\frac{x_{1}}{|x_{1}|}\cdot\nabla_{1}\psi_{1}\Big)
  e^{2(F-F_{1})}\psi_{1}
  \,dx_{2}\cdots dx_{N}.
\end{align}
Comparing
\eqref{eq:int_grad_psi}, \eqref{eq:not_alpha}, \eqref{eq:I3} and
\eqref{eq:I_5}, it can be seen that all the terms in \eqref{eq:except} 
stem from the function \(J_{1}\), namely from \(I_{1}\), \(I_{3}\),
and \(I_{5}\). 

All other terms in the decomposition of \(h\) 
are in \(C_{{\rm loc}}^{\alpha}(\R^{3})\). When integrating them
over \(\sphere^{2}\), we get
something continuous in \([0,\infty)\), according to
Proposition~\ref{prop:cont_int} above.

For the terms in \eqref{eq:except} we note that they are all of the
form
\begin{align}
  \label{eq:pull_out_omega}
   \int\frac{x_{1}}{|x_{1}|}\cdot{\bf K}(x_{1},x')
  \,dx'
  =\frac{x_{1}}{|x_{1}|}\cdot
  \int{\bf K}(x_{1},x')\,dx'.
\end{align}
In each case, we have
\begin{align}
   \label{eq:three_Holder}
    {\bf L}(x_1)=\big(L_{1}(x_{1}),L_{2}(x_{1}),L_{3}(x_{1})\big)=
  \int{\bf K}(x_{1},x')\,dx', L_{j}
  \in C^{\alpha}_{\text{{\rm
          loc}}}(\R^{3}). 
\end{align}
To see this, apply Lemma~\ref{lem:one} to each of the coordinate functions
\(L_{j}\), \(j=1,2,3\). The integrands are easily seen to satisfy
\eqref{eq:cond1} in each case, by the previous ideas. 
To get continuity in \([0,\infty)\) of the functions in  \eqref{eq:except}
we use \eqref{eq:pull_out_omega} and \eqref{eq:three_Holder}, and the 
following lemma:
\begin{prop}
  \label{prop:cont_int2}
  Assume \({\bf f}=(f_{1},f_{2},f_{3}), 
  f_{j}\in C^{\alpha}_{\text{{\rm loc}}}(\R^{3})\). 
 Then \(\bar f\in C^{0}([0,\infty))\), where
  \begin{align}
    \bar f(r)=\int_{\sphere^{2}}\big(\omega\cdot {\bf f}(r\omega)\big)
    \,d\omega.
    \nonumber
  \end{align}
\end{prop}
\begin{proof}
  The same as for Proposition~\ref{prop:cont_int}, noting that for all
  \(r\in[0,\infty)\) and fixed
  \(\omega\in\sphere^{2}\):
  \begin{align}
    \lim_{s\to r}\omega\cdot{\bf f}(s\omega)=\omega\cdot{\bf f}(r\omega).
    \nonumber
  \end{align}
\end{proof}
This holds even in the case \(r=0\), for which
\begin{align}
  \lim_{s\downarrow0}\int_{\sphere^{2}}\omega\cdot{\bf f}(s\omega)\,d\omega
  =\int_{\sphere^{2}}\omega\cdot{\bf f}(0)\,d\omega=0.
  \nonumber
\end{align}
  This proves that the functions in \eqref{eq:except} are in 
  \(C^{0}([0,\infty))\). Therefore
  \(\widetilde h\in C^{0}([0,\infty))\), which  
  finishes the proof of Proposition~\ref{prop:reg}.
\end{proof}

\textsl{Acknowledgement}: The authors wish to thank 
G.~Friesecke and S.~Fournais for stimulating discussions.

\providecommand{\bysame}{\leavevmode\hbox to3em{\hrulefill}\thinspace}


\end{document}